\RequirePackage{fix-cm}
\documentclass[smallextended]{svjour3}       
\smartqed  
\usepackage{amsfonts,color,graphicx}
 \journalname{-}
\begin{document}

\title{On the solution of a class of fuzzy system of linear equations 
}

\titlerunning{On the solution of a class of fuzzy system of linear equations}        

\author{Davod Khojasteh Salkuyeh
}

\authorrunning{D. K. Salkuyeh} 

\institute{D. K. Salkuyeh \at
              Faculty of Mathematical Sciences, University of Guilan, Rasht, Iran\\
              Tel.: +98-131-3233901\\
              Fax:  +98-131-3233509\\
              \email{khojasteh@guilan.ac.ir, salkuyeh@gmail.com}           
}

\date{Received: date / Accepted: date}

\maketitle

\begin{abstract}
In this paper, we consider the system of linear equations $Ax=b$, where $A\in \Bbb{R}^{n \times n}$ is a crisp H-matrix and $b$ is a fuzzy $n$-vector. We then investigate the existence and uniqueness of a fuzzy solution to this system. The results can also be used for the class of M-matrices and strictly diagonally dominant matrices. Finally, some numerical examples are given to illustrate the presented theoretical results.

\end{abstract}

\section{Introduction}
\label{SEC1}

Systems of linear equations play a major role in various field of science and engineering. In many applications, some of the system parameters are represented by fuzzy numbers rather than crisp numbers. Therefore, developing mathematical models and numerical procedures for the fuzzy system of linear equations (hereafter denoted by FSLE) would be of interest. In (Friedman et al. (1998)) the authors considered a general model for solving a FSLE whose coefficient matrix is crisp and its right-hand side is an arbitrary fuzzy vector. They stated some conditions for the existence of a unique fuzzy solution to FSLE by using  the embedding method (Friedman et al. (1998)) and converting the original system to a crisp linear system of equations. Later, several authors studied FSLE. In (Allahviranloo (2004)), the author used the Jacobi and Gauss-Seidel iterative methods to compute an approximate solution to Eq. (\ref{FSLE}). He also used the successive overrelaxation iterative method for solving FSLE in (Allahviranloo (2005)). Dehghan and Hashemi (2006) investigated the existence of a solution to Eq. (\ref{FSLE}) provided that the matrix $A$ is strictly diagonally dominant matrix with positive diagonal entries and then applied several iterative methods for solving FSLE. In (Hashemi et al. (2008)) the authors have studied FSLE under the condition that the coefficient matrix is an crisp M-matrix. Some block iterative methods to solve FSLE have been presented by Wang and Zheng (2007).

In this paper, we consider the FSLE of the the form
\begin{equation}\label{FSLE}
Ax=b,
\end{equation}
where $A\in \Bbb{R}^{n \times n}$ is a crisp H-matrix and $b$ is a fuzzy vector. We then investigate the existence and uniqueness of a solution to FSLE. To do so, we first review some notations, definitions and results that will be used in the next sections.

The vector $x=(x_1,x_2,\ldots,x_n)^T\in\mathbb{R}^n$ is said to be positive (nonnegative), if $x_i>0$ $(x_i\geq 0)$, $i=1,2,\ldots,n$. In this case, we write $x>0$ $(x\geq 0)$. Similar definitions can be written for matrices. A \textit{permutation matrix} is a matrix that has exactly one $1$ in each row or column and $0$s elsewhere. It is easy to show that if $P$ is a permutation matrix, then $P$ is nonsingular and $P^{-1}=P^T$.

\begin{definition}
(Axelsson (1996)) A matrix $A=(a_{ij})\in \Bbb{R}^{n\times n}$ is said to be an
M-matrix if $a_{ii}>0$ for $i=1,\ldots,n$, $a_{ij}\leq
0$, for $i\neq j$, $A$ is nonsingular and $A^{-1}\geq 0$. The \textit{comparison matrix} of $A$, denoted by $\mathcal{M}(A)=(m_{ij})\in \mathbb{R}^{n \times n}$, is defined by
\[
m_{ij}=\left\{
  \begin{array}{ll}
    ~~|a_{ii}|, & \textrm{if}~ i=j,\\[2mm]
    -|a_{ij}|, & \textrm{if}~ i\neq j.
  \end{array}
\right.
\]
$A$ is called an H-matrix if $\mathcal{M}(A)$ is an M-matrix.
\end{definition}

\begin{definition}\label{HmatrixDef}
(Axelsson (1996)) A matrix $A=(a_{ij})\in \Bbb{R}^{n\times n}$ is said to be
\textit{generalized strictly diagonally dominant} if
\begin{equation}\label{GSDDE}
|a_{ii}|x_i > \sum_{\stackrel{{j=1}}{j\neq i}}^n  |a_{ij}|x_j,\quad i=1,2,\ldots,n,
\end{equation}
for some positive vector $x=(x_1,x_2,\ldots,x_n)^T$. $A$ is called strictly diagonally dominant if Eq. (\ref{GSDDE}) is valid for $x=(1,1,\ldots,1)^T$.
\end{definition}

\begin{theorem}\label{HGSDDthm}
(Axelsson (1996), Lemma 6.4) A matrix $A\in \Bbb{R}^{n \times n}$ is an H-matrix if and only if $A$ is generalized strictly diagonally dominant.
\end{theorem}

Clearly, every H-matrix is nonsingular and every M-matrix is an H-matrix. From Theorem  \ref{HmatrixDef}, it is easy to see that every strictly diagonally dominant matrix is an H-matrix.

\indent This paper is organized as follows. In Section \ref{SEC2}, some basic definitions and results on FSLE are given. Our main results are drawn in Section \ref{SEC3}. Section \ref{SEC4} is devoted to some numerical examples illustrating the theoretical results. Some concluding remarks are given in Section \ref{SEC5}.

\section{Preliminaries} \label{SEC2}

In this section, we review some of the basic notations of fuzzy number arithmetic and fuzzy system of linear equations.

Following (Friedman et al. (1998)), a fuzzy number in parametric form is defined an ordered pair of functions $(\underline{{u}}(r),\overline{u}(r))$, $0\leq r \leq 1$, which satisfies the following requirements
\begin{enumerate}
  \item [(a)] $\underline{{u}}(r)$ is a bounded left continuous nondecreasing function over $[0,1]$;
  \item [(b)] $\overline{u}(r)$ is a bounded right continuous nonincreasing function over $[0,1]$;
  \item [(c)] $\underline{{u}}(r)\leq\overline{u}(r)$, $0\leq r\leq 1$.
\end{enumerate}
A crisp number $\alpha$ is represented by $\underline{{u}}(r)=\overline{u}(r)=\alpha$, $0\leq r\leq 1$.
A popular representation for fuzzy number is the trapezoidal representation $u=(x_0,y_0,\alpha,\beta)$
with defuzzifier interval $[x_0,y_0]$, left fuzziness $\alpha$ and right fuzziness $\beta$ (Hashemi et al. (2008)). The membership function of this trapezoidal number is as follows:
\[
u(x)=\left\{
    \begin{array}{ll}
      \frac{1}{\alpha}(x-x_0+\alpha), & ~~x_0-\alpha\leq x\leq x_0, \\[2mm]
      1, &~~ x_0\leq x\leq y_0,\\ [2mm]
      \frac{1}{\beta}(y_0-x+\beta), &~~ y_0\leq x\leq y_0+\beta, \\[2mm]
      0, &~~ otherwise.
    \end{array}
  \right.
\]
The parametric form of the number is
$$\underline{u}(r)=x_0-\alpha+\alpha r, \quad \overline{u}(r)=y_0+\beta-\beta r.$$
When $x_0=y_0$, a trapezoidal fuzzy number is reduced to a triangular fuzzy number.

To introduce the FSLE and defining its solution, we recall the arithmetic operations of arbitrary fuzzy numbers  $x=(\underline{{x}}(r),\overline{x}(r))$
and $y=(\underline{{y}}(r),\overline{y}(r))$, $0\leq r\leq 1$ and real number $k$ as follows
\begin{enumerate}
  \item [(a)] $x=y$ if and only if $\underline{{x}}(r)=\underline{{y}}(r)$ and $\overline{x}(r)=\overline{y}(r)$;\\[-2mm]
  \item [(b)] $x+y=(\underline{{x}}(r)+\underline{{y}}(r),\overline{x}(r)+\overline{y}(r))$; \\[-2mm]
  \item [(c)] $kx=\left\{
                   \begin{array}{ll}
                     (k\underline{{x}}(r),k\overline{x}(r)), &~~ k\geq0, \\[2mm]
                     (k\overline{x}(r),k\underline{{x}}(r)), &~~ k<0.
                   \end{array}
                 \right.$
 \end{enumerate}

\begin{definition}
(Friedman et al. (1998)) A fuzzy number  vector $x=(x_1,x_2,\ldots,x_n)^T$ where
$x_i=(\underline{x_i}(r),\overline{x_i}(r))$, $0\leq r\leq 1$, $i=1,2,\ldots,n$, is called a solution to FSLE (\ref{FSLE})
if
\begin{equation}\label{FuzzySolDef}
\left\{
\begin{array}{cc}
 \displaystyle\underline{\sum_{j=1}^{n} a_{ij}x_j}=\sum_{j=1}^{n} \underline{a_{ij}x_j}=\underline{y_i},\\
                         & \qquad i=1,2,\ldots,n.    \\
 \displaystyle\overline{\sum_{j=1}^{n} a_{ij}x_j}=\sum_{j=1}^{n} \overline{a_{ij}x_j}=\overline{y_i},
\end{array}
\right.
\end{equation}
\end{definition}

It is easy to see that Eq. (\ref{FuzzySolDef}) is equivalent to the $m \times m$ ($m=2n$) crisp system of linear equations
\begin{equation}\label{SmmEq}
\pmatrix{ S_1 & S_2 \cr S_2 & S_1 } \pmatrix{\underline{X} \cr \overline{X}}=\pmatrix{\underline{Y} \cr \overline{Y}},
\end{equation}
where $S_1=(\alpha_{ij})\in \Bbb{R}^{m \times m}$  in which
\[
\alpha_{ij}=\left\{
  \begin{array}{cl}
  a_{ij} , &~~ \textrm{if}~~ a_{ij}>0,  \\
  0      , &~~ \textrm{otherwise},
  \end{array}
\right.
\]
and $S_2=A-S_1$. Moreover, in Eq. (\ref{SmmEq}) we have $\underline{Y}=(\underline{y}_1,\underline{y}_2,\ldots,\underline{y}_n)^T$, $\overline{Y}=(\overline{y}_1,\overline{y}_2,\ldots,\overline{y}_n)^T$, $\underline{X}=(\underline{x}_1,\underline{x}_2,\ldots,\underline{x}_n)^T$ and
$\overline{X}=(\overline{x}_1,\overline{x}_2,\ldots,\overline{x}_n)^T$.

\begin{theorem}\label{ExiGen}
(Friedman et al. (1998), Theorem 1)
The coefficient matrix in Eq. (\ref{SmmEq}) is nonsingular if and only if $A=S_1+S_2$ and $S_1-S_2$ are both nonsingular.
\end{theorem}

\begin{theorem}\label{ExiGenFuzzy}
(Friedman et al. (1998), Theorem 3)
The components of the unique solution $X$ of Eq. (\ref{SmmEq}) represent a solution fuzzy vector to the system (\ref{FSLE}) for arbitrary $Y$ if and only if $S^{-1}$ is nonnegative, i.e., $S^{-1}\geq 0$.
\end{theorem}

\section{Main results} \label{SEC3}

In this section we investigate the existence and uniqueness of a fuzzy solution to (\ref{FSLE}) when its coefficient matrix is an H-marix.
For convenience we denote Eq. (\ref{SmmEq}) by $SX=Y$. We now state and prove the main theorem of the paper.

\begin{theorem}\label{MainThm}
If the coefficient matrix of system (\ref{FSLE}) is an H-matrix, then there exists a permutation matrix $P$ such that $PS$ is also an H-matrix.
\begin{proof}\rm
Let $A$ be an H-matrix. We construct a permutation matrix $P=(p_{ij})$ such that $\tilde{S}=PS$ is an H-matrix. To do so we first set $P=I$, where $I$ is the identity matrix of dimension $m$, and then  modify some of its entries as follows
\begin{equation}\label{PermutPdef}
\left\{
  \begin{array}{l}
    \textrm{For}~ i=1,\ldots,n,~\textrm{if}~S_{ii}=0,~ \textrm{then}~P_{ii}:=0~\textrm{and}~P_{i,i+n}:=1; \\[2mm]
    \textrm{For}~ i=n+1,\ldots,m, ~\textrm{if}~ S_{ii}=0,~ \textrm{then}~ P_{ii}:=0~ \textrm{and}~P_{i,i-n}:=1.
  \end{array}
\right.
\end{equation}
We claim that $P$ is the desired permutation matrix. Since $A$ is an H-matrix, according to Theorem \ref{HGSDDthm} there exists a positive vector $x=(x_1,x_2,\ldots,x_n)^T\in \Bbb{R}^n$ such that
\[
|a_{ii}|x_i > \sum_{\stackrel{{j=1}}{j\neq i}}^n |a_{ij}|x_j,\quad i=1,2,\ldots,n.
\]
Let
\[
\tilde{x}=(x^T,x^T)^T=(\tilde{x}_1,\ldots,\tilde{x}_n,\tilde{x}_{n+1},\ldots,\tilde{x}_m)^T.
\]
Obviously, we have $\tilde{x}>0$. It is not difficult to show that
\[
\left\{
  \begin{array}{ll}
    |\tilde{s}_{ii}| \tilde{x_i}= |a_{ii}|  x_i >  \displaystyle\sum_{\stackrel{{j=1}}{ j\neq i}}^n |a_{ij}|x_j=\displaystyle\sum_{\stackrel{{j=1}}{j\neq i}}^m |\tilde{s}_{ij}|\tilde{x}_j,    & ~~ \textrm{for}~ i=1,\ldots,n; \\ [2mm]
    |\tilde{s}_{ii}| \tilde{x_i}= |a_{i-n,i-n}|  x_{i-n} >  \displaystyle\sum_{\stackrel{{j=1}}{j\neq i-n}}^n |a_{i-n,j}|x_j=\displaystyle\sum_{\stackrel{{j=1}}{j\neq i}}^m |\tilde{s}_{ij}|\tilde{x}_j,    & ~~ \textrm{for}~ i=n+1,\ldots,m, \\ [2mm]

   \end{array}
\right.
\]
where $\tilde{s}_{ij}$'s are the entries of $\tilde{S}$. This shows that $\tilde{S}$ is an H-matrix. \qquad $\Box$
\end{proof}
\end{theorem}

\begin{corollary}\label{CorOne}
{\bf (a)} If the coefficient matrix of system (\ref{FSLE}) is an H-matrix with positive diagonal entries, then  $S$ is also an H-matrix.\\
{\bf (b)} If the coefficient matrix of system (\ref{FSLE}) is an M-matrix, then  $S$ is also an M-matrix.
\begin{proof}\rm
{\bf (a)} By refereing to Eq. (\ref{PermutPdef}), we see that if $A$ is an H-matrix with positive diagonal entries then $P=I$, where $I$ is the identity matrix. Therefore, $PS=S$ and as  a result $PS$ is an H-matrix.\\
{\bf (b)} Diagonal entries of an M-matrix are positive. On the other hand, every M-matrix is an H-matrix. Therefore, according to part (a), $P=I$ and as a result $PS=S$.  Obviously, diagonal entries of $S$ are positive and its offdiagonal entries are nonpositive. Therefore, we conclude that the matrix  $S$ is an M-matrix. \qquad $\Box$
\end{proof}
\end{corollary}

\begin{corollary}\label{CoreTwo}
{\bf (a)} If the coefficient matrix of system (\ref{FSLE}) is  strictly diagonally dominant, then there is a permutation matrix $P$  such that  $PS$ is also strictly diagonally dominant.\\
{\bf (b)} If the coefficient matrix of system (\ref{FSLE}) is  strictly diagonally dominant with positive diagonal entries, then $S$ is also strictly diagonally dominant.
\begin{proof}\rm
This corollary is an immediate result of the second part of  Definition \ref{HmatrixDef} together with Theorem \ref{MainThm} and Corollary \ref{CorOne}.
 \qquad $\Box$
\end{proof}
\end{corollary}

Some important observations can be posed here. From the presented theoretical results we see that if the matrix $A$ has any properties of being  H-matrix, M-matrix or strictly diagonally dominant then there is a permutation matrix $P$ (in the case that the diagonal entries of $A$ are positive we have $P=I$) such that $PS$ has the same property as $A$ and in any of these cases we have $PS$ is nonsingular. Since $P$ is nonsingular we can deduce that the matrix $S$ is nonsingular and the system $SX=Y$ has the unique solution $X^*=S^{-1}Y$.

If $A$ is an H-matrix or strictly diagonally dominant matrix then there it does not guarantee for the computed solution to Eq. (\ref{FSLE}) being a fuzzy vector. However, when the coefficient matrix of Eq. (\ref{FSLE}) is an M-matrix, according to the second part of Corollary \ref{CorOne}, $S$ is an M-matrix and  we have $S^{-1}\geq 0$. Therefore, by Theorem \ref{ExiGenFuzzy} the solution of system (\ref{FSLE}) is a fuzzy vector. Moreover, as mentioned in (Salkuyeh (2011)), when $A$ is an M-matrix then $S$, can be partitioned as
\begin{equation}\label{MmatrixS}
S=\pmatrix{D    &   E \cr
           E    &   D },
\end{equation}
where  $D={\rm diag}(a_{11},a_{22},\ldots,a_{nn})$ and $E=A-D$. Therefore several iterative methods such as the AOR iterative method (Hadjidimos (1978), Salkuyeh (2011))), or direct methods such the LU factorization of $S$ can be exploited  to solve the system $SX=Y$. The book (Meurant (1999)) is a comprehensive reference for different methods to solve such systems.

The last comment which we present here is due to the study of a method to solve the more general case that
the coefficient matrix $A$ is an H-matrix. If $A$ is an H-matrix, then by Theorem \ref{MainThm} the matrix $\tilde{S}$ is an H-matrix. It is well-known that if $A$ is an H-matrix then its LU factorization exists (Meurant (1999), Theorem 2.17). Therefore, one can solve $\tilde{S}X=PY$ instead of $SX=Y$. In fact, since $\tilde{S}$ is an H-matrix, we can use its LU factorization to solve $\tilde{S}X=PY$. In this paper we will use this method for the presented examples.

\section{Numerical examples}\label{SEC4}

{In this section, some numerical experiments are given to illustrate the theoretical results presented in this paper. The triangular fuzzy numbers are used in all of the following numerical examples.

\begin{example}
In the first example, we consider the FSLE  $Ax=b$, where
\[
A=\pmatrix{~~2 &   -3  & ~~1 \cr
           -1  &   -3  &  -1 \cr
           -1  & ~~ 2  & ~~5},\qquad
b=\pmatrix{(  -11 + 9r  ,    7 - 9r  ) \cr
           (  -22 + 8r  ,   -6 - 8r  ) \cr
           (    9 +10r  ,   29 -10r  ) }.
\]
It is easy to see that $A$ is an H-matrix. The corresponding matrix $S$ and its associated permutation matrix $P$ are
\[
S=\pmatrix{~~ 2  &  ~~ 0  &   ~~1  & ~~0  &  -3  & ~~0 \cr
           ~~ 0  &  ~~ 0  &   ~~0  &  -1  &  -3  &  -1 \cr
           ~~ 0  &  ~~ 2  &   ~~5  &  -1  & ~~0  & ~~0 \cr
           ~~ 0  &    -3  &   ~~0  & ~~2  & ~~0  & ~~1 \cr
             -1  &    -3  &    -1  & ~~0  & ~~0  & ~~0 \cr
             -1  &  ~~ 0  &   ~~0  & ~~0  & ~~2  & ~~5 },\qquad
P=\pmatrix{1  &   0  &   0  &   0  &   0  &   0 \cr
           0  &   0  &   0  &   0  &   1  &   0 \cr
           0  &   0  &   1  &   0  &   0  &   0 \cr
           0  &   0  &   0  &   1  &   0  &   0 \cr
           0  &   1  &   0  &   0  &   0  &   0 \cr
           0  &   0  &   0  &   0  &   0  &   1}.
\]
We have
\[
\tilde{S}=PS= \pmatrix{  ~~2 &  ~~0  & ~~1  & ~~0  &  -3  & ~~0 \cr
                          -1 &   -3  &  -1  & ~~0  & ~~0  & ~~0 \cr
                         ~~0 &  ~~2  & ~~5  &  -1  & ~~0  & ~~0 \cr
                         ~~0 &   -3  & ~~0  & ~~2  & ~~0  & ~~1 \cr
                         ~~0 &  ~~0  & ~~0  &  -1  &  -3  &  -1 \cr
                          -1 &  ~~0  & ~~0  & ~~0  & ~~2  & ~~5 },
\]
which is an H-matrix. We see that Theorem \ref{MainThm} holds here. Therefore, we  solve the system $\tilde{S}X=PY$, by using the LU factorization of $\tilde{S}$. This results in the exact solution of the system as
\[
x^*=\pmatrix{(1+ r,3-r ) \cr
             (1+2r,5-2r) \cr
             (2+r,4-r)}.
\]
\end{example}


\begin{example}
In this example we consider the FSLE  $Ax=b$, where
\[
A=\pmatrix{  2 & -1  &  1 \cr
            -3 &  5  &  1 \cr
            -2 & -3  &  4 },\qquad
b=\pmatrix{(    -2 + 9r  ,   13 - 6r) \cr
           (    -7 +12r  ,   25 -20r ) \cr
           (   -23 +14r  ,    4 -13r ) }.
\]
Here the matrix $A$ is an H-matrix with positive diagonal entries. The corresponding matrix $S$ and its associated permutation matrix $P$ are
\[
S=\pmatrix{    ~~2  & ~~0  &   1  & ~~0  &  -1  &   0 \cr
               ~~0  & ~~5  &   1  &  -3  & ~~0  &   0 \cr
               ~~0  & ~~0  &   4  &  -2  &  -3  &   0 \cr
               ~~0  &  -1  &   0  & ~~2  & ~~0  &   1 \cr
                -3  & ~~0  &   0  & ~~0  & ~~5  &   1 \cr
                -2  &  -3  &   0  & ~~0  & ~~0  &   4 },\qquad
P=I_{6 \times 6}.
\]
According to the part (a) of Corollary \ref{CorOne} the matrix $S$ is an H-matrix. Therefore, we exploit its LU factorization for solving the system $SX=Y$.
In this case, we obtain
\[
x^*=\pmatrix{(1+3r  , 6- 2r ) \cr
             (2+r   , 5-2r  ) \cr
             (1+r   , 3- r  )}.
\]
\end{example}


\begin{example}
Our last example is devoted to the FSLE  $Ax=b$, where
\[
A=\pmatrix{~~6 & -1  & -1 \cr
            -1 & ~~2 & -1 \cr
            -1 & -1  & ~~1  },\qquad
b=\pmatrix{( -18 +16r    ,  8 -10r ) \cr
           (  -8 + 8r    ,  6 - 6r ) \cr
           (  -3 + 4r    ,  8 - 7r ) }.
\]
Here the coefficient matrix $A$ is an M-matrix. Since its diagonal entries are positive the corresponding permutation matrix is  $P=I_{6 \times 6 }$ and
\[
S=\pmatrix{    ~~6  & ~~0  & ~~0  & ~~0  &  -1  &  -1   \cr
               ~~0  & ~~2  & ~~0  &  -1  & ~~0  &  -1   \cr
               ~~0  & ~~0  & ~~1  &  -1  &  -1  & ~~0   \cr
               ~~0  &  -1  &  -1  & ~~6  & ~~0  & ~~0   \cr
                -1  & ~~0  &  -1  & ~~0  & ~~2  & ~~0   \cr
                -1  &  -1  & ~~0  & ~~0  & ~~0  & ~~1   },
\]
As seen, $S$ is of the form (\ref{MmatrixS}). According to the part (b) of Corollary \ref{CorOne} the matrix  $S$ is an
M-matrix and its LU factorization exists. Invoking the LU factorization of $S$ the exact solution of $Ax=b$ is obtained as
\[
x^*=\pmatrix{( -1+ 2r , 2- r  ) \cr
             (  1+ 2r , 4- r  ) \cr
             (  3+2r  , 8-3r )}.
\]
\end{example}

\section{Conclusion}\label{SEC5}

In this paper we have verified the existence and uniqueness of a solution to the fuzzy system of linear equations $Ax=b$ where $A$ in a crisp H-matrix and $b$ is a fuzzy vector. Our main result in special cases can be applicable for the class of M-matrices and strictly diagonally dominant matrices. Some numerical examples have been shown to illustrate the theoretical results.

\bigskip

\noindent \textbf{References}

\bigskip

\noindent Allahviranloo T 2004 Numerical methods for fuzzy system of linear equations, \emph{Applied}

\emph{Mathematics and Computation} 155: 493--502.\\[-3mm]

\noindent Allahviranloo T 2005 Successive overrelaxation iterative method for fuzzy system of linear

 equations, \emph{Applied  Mathematics and Computation} 162: 189--196. \\[-3mm]

\noindent Axelsson O 1996  \emph{Iterative solution methods}, Cambridge University Press, Cambridge, 1996.\\[-3mm]

\noindent Dehghan M, Hashemi B 2006 Iterative solution of fuzzy linear systems, \emph{Applied  Mathematics}

\emph{ and Computation} 175: 645--674.\\[-3mm]

\noindent Ezzati R 2011 Solving fuzzy linear systems, \emph{Soft Computing} 15: 193--197.\\[-3mm]

\noindent Friedman M,  Ming M, Kandel A 1998 Fuzzy linear systems, \emph{Fuzzy Sets and Systems}

 96: 201--209.\\[-3mm]

\noindent Hadjidimos A 1978 Accelerated overrelaxation method, \emph{Mathematics of Computation}

32: 149--157.\\[-3mm]

\noindent  Hashemi M S,   Mirnia M K,  Shahmorad S 2008, Solving  fuzzy linear systems by using the

Schur complement when coefficient matrix is an M-matrix, \emph{Iranian Journal of Fuzzy }

\emph{Systems} 5: 15--29.\\[-3mm]

\noindent Salkuyeh D K 2011 On the solution of the fuzzy Sylvester matrix equation,

\emph{Soft Computing} 15: 953--961.\\[-3mm]

\noindent Meurant G 1999  \emph{Computer solution of large linear systems},  North-Holland, Amsterdam.\\[-3mm]

\noindent Wang K, Zheng B 2007 Block iterative methods for fuzzy linear systems, \emph{Journal of Applied}

\emph{Mathematics and Computing} 25: 119--136.\\[-3mm]

\noindent Zadeh L A 1965 Fuzzy sets, \emph{Information and Control} 8: 338--353.

\end{document}